\documentclass[12pt,a4paper]{article}

\usepackage[T1]{fontenc}
\usepackage[utf8]{inputenc}
\usepackage{lmodern}
\usepackage[margin=1in]{geometry}
\usepackage{amsmath,amssymb,amsthm}
\usepackage{mathtools}
\usepackage{graphicx}
\usepackage{booktabs}
\usepackage{enumitem}
\usepackage[colorlinks=true,linkcolor=blue,citecolor=blue,urlcolor=blue]{hyperref}

\theoremstyle{plain}
\newtheorem{theorem}{Theorem}[section]
\newtheorem{lemma}[theorem]{Lemma}

\newtheorem{corollary}[theorem]{Corollary}

\theoremstyle{definition}

\theoremstyle{remark}

\newcommand{\cM}{\mathcal{M}}
\newcommand{\cG}{\mathcal{G}}
\newcommand{\cD}{\mathcal{D}}
\newcommand{\p}{\mathbf{p}}

\newcommand{\bbp}{\mathbb{P}}

\newcommand{\Poi}{\operatorname{Poi}}

\newcommand{\PP}{\bbp}
\DeclareMathOperator{\rk}{rk}

\title{Asymptotic enumeration of minimally $d$-rigid graphs}
\author{Motty Gruda \and Yuval Peled\thanks{M.G. and Y.P. were partially supported by the Israel Science Foundation grant ISF-3464/24.}}
\date{}

\begin{document}

\maketitle

\begin{abstract}
For every fixed $d\ge2$, we determine an explicit constant $\kappa_d>0$
such that the number of minimally $d$-rigid graphs on $n$ labelled
vertices is $n^{dn}(\kappa_d+o(1))^n$ as $n\to\infty$.
Our proof relies on separately enumerating the $(d+1)$-core of the graph and the
remaining edges, which we call the {mantle}.

To enumerate rigid cores, we show that a graph chosen uniformly from
$n$-vertex graphs with minimum degree at least $d+1$ and exactly
$dn-\binom{d+1}{2}$ edges is $d$-rigid with probability $(1-o(1))^n$.
To enumerate mantles, we use a bijection with generalized parking
functions. 

The enumeration also yields structural results for a
uniformly random minimally $d$-rigid graph. In particular, the fraction
of vertices in its $(d+1)$-core is exponentially concentrated around
$\PP(X\ge d+1\mid X\ge d)$, where $X$ is Poisson distributed with parameter chosen
so that $\mathbb E[X\mid X\ge d+1]=2d$.
\end{abstract}

\section{Introduction}
\label{sec:introduction}
Generic graph rigidity is a well-studied topic in discrete geometry. A $d$-framework is a pair $(G,\p)$ consisting of a graph $G=(V,E)$ and a realization $\p:V\to\mathbb{R}^d$ of its vertices. The realization is called generic if its coordinates are algebraically independent over $\mathbb{Q}$. 
A graph $G=(V,E)$ is called \emph{$d$-rigid} if, for a generic realization
$\p:V\to\mathbb{R}^d$, every sufficiently small motion  of the vertices starting from $\p$ that preserves all edge lengths is induced by a Euclidean isometry.
A graph is \emph{minimally $d$-rigid} if it is $d$-rigid and deleting any edge permits a nontrivial motion that preserves the lengths of all remaining edges. 
Going back to the physicist Maxwell~\cite{Maxwell1864}, it is known that every minimally $d$-rigid graph on $n\ge d$ vertices has exactly $dn-\binom{d+1}{2}$ edges~\cite{asimow1}.
Minimally rigid graphs are central to rigidity theory. For example, the celebrated Laman's theorem~\cite{Laman1970OnGA}
characterizes minimally $2$-rigid graphs by a sparsity condition.
Finding a comparably simple combinatorial characterization in dimensions $d\ge3$ remains a major open problem. 

For $d=1$, the minimally rigid graphs are spanning trees, and Cayley's formula shows that there are  $n^{n-2}$ spanning trees on $n$ labelled vertices. In this paper, we find the asymptotic number of minimally $d$-rigid graphs on $n$ labelled vertices for every fixed $d\ge2$ as $n\to\infty$, and study their typical structure.

Let $\cM_{n,d}$ denote the family of minimally $d$-rigid graphs
on the vertex set $[n]=\{1,\ldots,n\}$. It is not hard to see that $|\cM_{n,d}|=(\Theta(n))^{dn}$. 
Indeed, an elementary upper bound  $$|\cM_{n,d}|\le \binom{\binom n2}{dn-\binom{d+1}{2}}$$
  is given by the number of $n$-vertex $\left(dn-\binom{d+1}{2}\right)$-edge graphs. For a lower bound, we use the fact that if $G$ is minimally $d$-rigid, then adding a new vertex and joining it to $d$ distinct vertices of $G$ produces a larger minimally $d$-rigid graph. This operation is called a \emph{$0$-extension}~\cite{Henneberg1911}. Let $\cM_{n,d}(0)$ denote the family of graphs on $[n]$ obtainable from $K_d$ by $0$-extensions. To find a lower bound for $|\cM_{n,d}(0)|$, we
start from a fixed labelled copy of $K_d$, order the remaining $n-d$
vertices, and join each new vertex to its immediate predecessor and to
$d-1$ other earlier vertices. The construction is injective because repeatedly deleting
the unique degree-$d$ vertex outside the fixed initial $K_d$ recovers the order, and thus
\[
|\cM_{n,d}(0)| \ge (n-d)!\cdot\prod_{i=1}^{n-d} \binom{d+i-2}{d-1}\,.
\]
A straightforward calculation using Stirling's formula therefore gives
\begin{equation}
    \left(\frac{1}{e^d(d-1)!}+o(1)\right)^n n^{dn}
    \le |\cM_{n,d}|
    \le \left(\left(\frac{e}{2d}\right)^d+o(1)\right)^n n^{dn}.
    \label{eq:elementary-bounds}
\end{equation}
Due to the fact that each new vertex is attached to its immediate predecessor, the
lower bound in \eqref{eq:elementary-bounds} counts only a subfamily of $\cM_{n,d}(0)$. Angel and Kolesnik's work on minimally susceptible graphs in $d$-neighbour bootstrap percolation~\cite{Angel2016ThresholdsFC} gives a precise asymptotic count:
\begin{equation}
    |\cM_{n,d}(0)| =
    \left(\frac{1}{e^{d-1}(d-1)!}+o(1)\right)^n n^{dn}.
    \label{eq:zero-extension-count}
\end{equation}
The lower and upper bounds above share the factor $n^{dn}$ and differ only in the constant raised to the power $n$. Our main result determines this constant for $|\cM_{n,d}|$.

Let $\lambda_d^*$
be the unique positive $\lambda$ for which
\begin{equation}
    \mathbb{E}_{X\sim\Poi(\lambda)}[X\mid X\ge d+1]
    =2d.
    \label{eq:poisson-parameter}
\end{equation}
Define
\begin{equation}
    \kappa_d=
    \frac{e^{-\lambda_d^*/2}}
         {d!\left(1-\lambda_d^*/(2d)\right)}.
    \label{eq:growth-constant}
\end{equation}

\begin{theorem}\label{thm:enumeration}
For every fixed $d\ge 2$, as $n\to\infty$,
\[
    |\cM_{n,d}|=(\kappa_d+o(1))^n n^{dn}.
\]
\end{theorem}

The enumerative result in Theorem \ref{thm:enumeration} relies on a canonical decomposition of (almost all) minimally $d$-rigid graphs $G$ into two independent parts:
\begin{enumerate}
\item The \emph{$(d+1)$-core}: the unique maximal subgraph of $G$ with minimum degree at least $d+1$.
\item The \emph{mantle}: the set of edges outside the $(d+1)$-core admits a unique acyclic orientation in which every vertex outside the core has out-degree $d$, while every vertex in the core has out-degree $0$.
\end{enumerate}

Note that the $(d+1)$-core may be empty, in which case $G\in\cM_{n,d}(0)$, and such graphs are enumerated by \eqref{eq:zero-extension-count}.
Otherwise, the core is itself minimally $d$-rigid, and $G$ can be obtained from it by a sequence of $0$-extensions. 
Once the core vertex set is fixed, the core and the mantle can be
chosen independently. Attaching any such mantle to a minimally $d$-rigid
graph on that set with minimum degree at least $d+1$ produces a minimally
$d$-rigid graph.
Consequently,
\begin{equation}
|\cM_{n,d}|
=
|\cM_{n,d}(0)|
+
\sum_{k=d+2}^{n}
\binom{n}{k}
|\cM_{k,d}^{(d+1)}|
|\cD_{k,n-k}^{(0,d)}|\,,
\label{eq:rigid-core-decomposition}
\end{equation}
where $\cM_{k,d}^{(d+1)}$ denotes the
family of minimally $d$-rigid graphs on $[k]$ with minimum degree
at least $d+1$, and $\cD_{k,n-k}^{(0,d)}$ denotes the family of
directed acyclic graphs on $[n]$ in which vertices $1,\ldots,k$
have out-degree $0$ and the remaining vertices have out-degree $d$.

The asymptotic enumeration of graphs with a prescribed degree sequence is a well-developed topic~\cite{LiebenauWormald2024}.
Let $\cG^{(d+1)}_{n}$ denote the family of  graphs on $[n]$ with minimum degree at least $d+1$ and exactly $dn-\binom{d+1}{2}$ edges. Pittel and Wormald provide an asymptotic enumeration of this family given by
$$
|\cG^{(d+1)}_{n}|=(b_d+o(1))^n n^{dn}\,,
$$
where
\[
    b_d=\frac{(2d/e)^d}
               {d!\,(\lambda_d^*)^{d-1}(2d-\lambda_d^*)}\,.
\]
They also characterize its typical degree sequence in terms of a truncated Poisson distribution~\cite{kcore_enum}.
Building on recent results by the second author on rigidity of sparse random graphs~\cite{Peled2026}, we show that imposing the additional requirement of $d$-rigidity reduces the count by at most a subexponential factor.
\begin{theorem}\label{thm:rigid-cores}
Fix $d\ge 2$. Then, as $n\to\infty$,
$|\cM_{n,d}^{(d+1)}|=(b_d+o(1))^n n^{dn}\,.$
\end{theorem}
 In other words, if $G$ is chosen uniformly from $\cG^{(d+1)}_{n}$ then the probability that $G$ is minimally $d$-rigid is $(1-o(1))^n$.
We conjecture that this probability
is in fact $1-o(1)$, but our current methods are insufficient to prove this.

The following asymptotic enumeration of the mantle complements Theorem~\ref{thm:rigid-cores}.

\begin{theorem}
\label{thm:acyclic-attachments}
Fix $d\ge2$. For integers $d\le k\le n$, put $\alpha=k/n$.
Then, uniformly in $k$ as $n\to\infty$,
\[
    |\cD_{k,n-k}^{(0,d)}|=\binom nd^{n-k}\exp\bigl(-nJ_d(\alpha)+o(n)\bigr),
\]
where
\[
    J_d(\alpha)=
    \begin{cases}
        (d-1)-d\alpha
        +(d-1)\alpha\log\dfrac{d\alpha}{d-1}
        -(1-\alpha)\log\bigl(d(1-\alpha)\bigr),
        & 0\le\alpha\le 1-\dfrac1d,\\[6pt]
        0, & 1-\dfrac1d\le\alpha\le1.
    \end{cases}
\]
\end{theorem}

The main idea in the proof is a bijection between $\cD_{k,n-k}^{(0,d)}$ and
generalized parking functions~\cite{KungYan2003}. These are sequences $(f_1,...,f_{n-k})$ of positive integers such that for every $i=1,\ldots,n-k$, the sequence takes values that are at most $\binom{k+i-1}{d}$ at least $i$ times.

A direct analytical argument combining \eqref{eq:rigid-core-decomposition} and Theorems~\ref{thm:rigid-cores}
and~\ref{thm:acyclic-attachments} gives 
Theorem~\ref{thm:enumeration} as well as an exponential concentration of
the core size of a random minimally $d$-rigid graph around $\alpha_d n$, where
\begin{equation}
    \alpha_d
    =\PP_{X\sim\Poi(\lambda_d^*)}(X\ge d+1\mid X\ge d)
    =\frac{\lambda_d^*}{2d}.
    \label{eq:core-fraction}
\end{equation}

\begin{corollary}\label{cor:rigid-core-concentration}
Let $G$ be chosen uniformly from $\cM_{n,d}$ and let $C_{d+1}(G)$
denote its $(d+1)$-core.
For every $\varepsilon>0$, there exists $\delta=\delta(d,\varepsilon)>0$
such that, for all sufficiently large $n$,
\[
    \PP\!\left(
        \left|\,|V(C_{d+1}(G))|-\alpha_d n\,\right|>\varepsilon n
    \right)\le e^{-\delta n}.
\]
\end{corollary}

We also prove exponential concentration of the $0$-extension level sizes of $G$ in
Theorem~\ref{thm:level-concentration} below. Here level-$0$ is the $(d+1)$-core of $G$
and each subsequent level consists of the remaining vertices with
at least $d$ neighbors in earlier levels.

The remainder of the paper is organized as follows.
Section~\ref{sec:rigid-cores} proves Theorem~\ref{thm:rigid-cores}
on the abundance of minimally $d$-rigid cores.
Section~\ref{sec:attachments} proves Theorem~\ref{thm:acyclic-attachments}
by counting mantles through generalized parking functions.
Section~\ref{sec:main-proofs} proves Theorem~\ref{thm:enumeration}
and Corollary~\ref{cor:rigid-core-concentration} and discusses further
properties of typical minimally $d$-rigid graphs.

\section{Many minimally rigid cores}
\label{sec:rigid-cores}
\label{alt:rigid-cores}

We first obtain many $d$-independent cores from subcritical random
graphs and complete them by adding edges, with a controlled exponential
loss. This gives the desired count along a dense subsequence of vertex
counts; counting $1$-extensions then extends the limit to every $n$.

\subsection{Infinitesimal rigidity and the rigidity matroid}

For a graph $G=(V,E)$ and a generic placement
$p:V\to\mathbb{R}^d$, the rigidity matrix $R(G,p)$ has one row for each
edge and $d$ columns for each vertex. The row corresponding to $uv\in E$
has $p(u)-p(v)$ in the columns of $u$, $p(v)-p(u)$ in the columns of $v$,
and zeros elsewhere. Its kernel consists of the infinitesimal motions
of $(G,p)$, the velocity assignments $\xi:V\to\mathbb{R}^d$ satisfying
\[
    (p(u)-p(v))\cdot(\xi(u)-\xi(v))=0
    \qquad\text{for every }uv\in E.
\]
A framework is infinitesimally rigid if all such motions are
infinitesimal Euclidean motions. For generic $p$, rigidity is
equivalent to infinitesimal rigidity~\cite{asimow1}.
The rank of $R(G,p)$ is independent of the generic placement, and
we denote it by $\rk_d(G)$. For $|V|\ge d+1$,
\[
    \rk_d(G)\le\min\left\{|E|,\,d|V|-\binom{d+1}{2}\right\},
\]
and $G$ is $d$-rigid if and only if equality holds in the second bound.

The rows indexed by $E(K_n)$ define the \emph{generic $d$-rigidity
matroid}. An edge set is \emph{$d$-independent} if its rank equals its
cardinality. Its bases are precisely the minimally $d$-rigid graphs
on $[n]$, each with $dn-\binom{d+1}{2}$ edges.
The \emph{closure} of a graph $G$ on $[n]$ is
\[
    \operatorname{cl}_d(G)
    :=\{e\in E(K_n):\rk_d(G+e)=\rk_d(G)\}\supseteq E(G).
\]

\begin{lemma}[Closure bound]
\label{alt:closure-bound}
For every graph $G$ on $n\ge d+1$ vertices,
\[
    \frac{|\operatorname{cl}_d(G)|}{\binom n2}
    \le\frac{\rk_d(G)}{dn-\binom{d+1}{2}}.
\]
\end{lemma}

\begin{proof}
Let $B$ be a fixed minimally $d$-rigid graph and let $B_\sigma$ be its image under a
uniformly random permutation of $[n]$. On the one hand, the edge subset
$B_\sigma\cap\operatorname{cl}_d(G)\subseteq\operatorname{cl}_d(G)$ is $d$-independent and thus has at most $\rk_d(G)$ edges.
On the other hand, each edge belongs to $B_\sigma$ with probability $|B|/\binom n2$. Taking expectations proves the inequality.
\end{proof}

\subsection{Random graph preliminaries}

Throughout this section, $d\ge2$ is fixed and
$m_n=dn-\binom{d+1}{2}$. For variable edge counts, let
$\cG^{(d+1)}_{n,m}$ be the family of graphs on $[n]$ with $m$ edges
and minimum degree at least $d+1$. Let $\mathcal I^{(d+1)}_{n,m}$
be its subfamily consisting of $d$-independent graphs.
Thus $\cG^{(d+1)}_n=\cG^{(d+1)}_{n,m_n}$, in the notation of the
introduction.

We write $G(N,p)$ for the binomial random graph on $[N]$, and
$\pi_j(t)=\PP(\Poi(t)\ge j)$. Define
\[
    \gamma_d=\min_{t>0}\frac{t}{\pi_d(t)},
    \qquad
    c_d=\frac{\lambda_d^*}{\pi_d(\lambda_d^*)}.
\]
Here $\gamma_d$ is the emergence threshold of the $(d+1)$-core,
while $c_d$ is the rigidity-independence threshold.
For $c>\gamma_d$, let $\lambda(c)$ be the largest positive solution
of $\lambda=c\pi_d(\lambda)$.

\begin{theorem}[{\cite{kcore_pittel,JansonSimplekcore}}]
\label{alt:random-core}
Let $G_N\sim G(N,c/N)$, where $c>0$ is fixed, and put
$C_N=C_{d+1}(G_N)$.
If $c<\gamma_d$, then $C_N$ is empty a.a.s.
If $c>\gamma_d$, then, in probability,
\[
    \frac{|V(C_N)|}{N}\longrightarrow\pi_{d+1}(\lambda(c)),
    \qquad
    \frac{|E(C_N)|}{N}\longrightarrow
       \frac{\lambda(c)\pi_d(\lambda(c))}{2}.
\]
Conditional on its vertex set and its number of edges, $C_N$ is
uniform among graphs on that set with the same edge count and
minimum degree at least $d+1$.
\end{theorem}

Consequently, for $c>\gamma_d$, the limiting average degree of the
core is
\[
    \rho(c)=\frac{\lambda(c)\pi_d(\lambda(c))}
                   {\pi_{d+1}(\lambda(c))}
    =\mathbb E[\Poi(\lambda(c))\mid\Poi(\lambda(c))\ge d+1].
\]
This is continuous and strictly increasing in $c$. Moreover,
$\gamma_d<c_d$, $\lambda(c_d)=\lambda_d^*$, and $\rho(c_d)=2d$.
Thus the threshold $c_d$ is exactly where the nonempty core changes
from having fewer than $d$ edges per vertex to having more than $d$:
a.a.s.,
\[
    |E(C_N)|<d|V(C_N)|\quad\text{if }\gamma_d<c<c_d,
    \qquad
    |E(C_N)|>d|V(C_N)|\quad\text{if }c>c_d.
\]
These are standard properties of the core and orientability
thresholds~\cite{kcore_pittel,fernholz2007k}.

The following theorem of the second author implies that the
$(d+1)$-core of a subcritical $G(N,c/N)$ is a.a.s.\ $d$-independent,
a fact crucial to our proof. For $d=2$, this was already proved by
Kasiviswanathan, Moore and Theran~\cite{rigid_giant_component}.

\begin{theorem}[{\cite[Theorem~1.1(1)]{Peled2026}}]
\label{alt:independence-threshold}
For every fixed $0<c<c_d$, the random graph
$G_N\sim G(N,c/N)$ is a.a.s.\ $d$-independent.
\end{theorem}

We also need Pittel and Wormald's enumeration at average degrees
other than $2d$. For $\rho>d+1$, let $\lambda_\rho$ be the unique
positive solution of
\[
    \mathbb E[\Poi(\lambda_\rho)\mid\Poi(\lambda_\rho)\ge d+1]
    =\frac{\lambda_\rho\pi_d(\lambda_\rho)}
           {\pi_{d+1}(\lambda_\rho)}=\rho,
\]
and extend the notation $b_d$ from the introduction by putting
\[
    b_d(\rho)=\left(\frac{\rho}{e}\right)^{\rho/2}
       \frac{e^{\lambda_\rho}\pi_{d+1}(\lambda_\rho)}
            {\lambda_\rho^{\rho}}.
\]
The conditional mean above is continuous and strictly increasing
from $d+1$ to infinity, so $\lambda_\rho$ and $b_d(\rho)$ are
continuous. At $\rho=2d$, we have $\lambda_{2d}=\lambda_d^*$ and
$b_d(2d)=b_d$.

\begin{theorem}[Pittel and Wormald~\cite{kcore_enum}, Theorem~2(a)]
\label{alt:core-enumeration}
If $n\to\infty$ and $2m/n\to\rho>d+1$, then
\[
    |\cG^{(d+1)}_{n,m}|=(b_d(\rho)+o(1))^n n^m.
\]
\end{theorem}

For the remainder of the proof, write
\[
    a_n=\frac1n\log\frac{|\cM^{(d+1)}_{n,d}|}{n^{dn}},
\]
with $a_n=-\infty$ when the family is empty.
Theorem~\ref{thm:rigid-cores} asserts that $a_n\to\log b_d$.
The inclusion $\cM^{(d+1)}_{n,d}\subseteq\cG^{(d+1)}_n$ and the
preceding enumeration formula already give
\begin{equation}
    \limsup_{n\to\infty}a_n\le\log b_d.
    \label{alt:core-upper-bound}
\end{equation}

\subsection{Convergence along a dense subsequence}

A strictly increasing sequence $(n_i)$ is \emph{dense} if
$n_{i+1}/n_i\to1$. We also call its set of terms dense.
Equivalently, a set $A\subseteq\mathbb N$ is dense if, for every
$\eta>0$ and all sufficiently large $x$,
\[
    A\cap(x,(1+\eta)x]\ne\varnothing.
\]
Indeed, for consecutive terms $u_i\le x<u_{i+1}$, density puts
$u_{i+1}$ in this interval; conversely, apply this property at $x=u_i$.

We will use the following diagonalization twice: first to transfer
independence to a uniform core, and then to let its density approach $2d$.

\begin{lemma}[Diagonalization over dense subsequences]
\label{alt:dense-diagonal}
Let $A_1,A_2,\ldots$ be dense sets of positive integers, and let
$R_1,R_2,\ldots$ be arbitrary thresholds. There are a dense sequence
$(n_i)$ and nondecreasing integers $t_i\to\infty$ such that
$n_i\in A_{t_i}$ and $n_i>R_{t_i}$ for every $i$.
\end{lemma}

\begin{proof}
Choose strictly increasing integers $T_t>R_t$ so large that
\[
    A_t\cap(x,(1+1/t)x]\ne\varnothing
    \qquad\text{for every }x\ge T_t.
\]
Start with $t_1=1$ and any $n_1\in A_1$ with $n_1\ge T_1$.
Given $n_i\ge T_{t_i}$, set $t_{i+1}=t_i+1$ if
$n_i\ge T_{t_i+1}$, and otherwise set $t_{i+1}=t_i$.
Choose the least element $n_{i+1}$ of $A_{t_{i+1}}$ exceeding $n_i$.
Then
\[
    n_i<n_{i+1}\le(1+1/t_{i+1})n_i,
    \qquad n_{i+1}\ge T_{t_{i+1}}.
\]
The sequence $(n_i)$ is strictly increasing, hence tends to infinity.
Consequently $t_i$ cannot remain bounded: if it eventually stayed at
$t$, then $n_i\ge T_{t+1}$ would eventually force an increase.
Thus $t_i\to\infty$ and $n_{i+1}/n_i\to1$.
\end{proof}

\begin{lemma}[Transfer to the uniform model]
\label{alt:core-transfer}
Fix $c>\gamma_d$, and let $\mathcal P$ be a graph property invariant
under relabelling. Suppose that the $(d+1)$-core of
$G_N\sim G(N,c/N)$ has $\mathcal P$ a.a.s. Then there are integers
$(n_i,m_i)$ such that
\[
    \frac{n_{i+1}}{n_i}\longrightarrow1,
    \qquad
    \frac{2m_i}{n_i}\longrightarrow
       \rho(c),
\]
and a uniform graph in $\cG^{(d+1)}_{n_i,m_i}$ has $\mathcal P$
with probability tending to $1$.
\end{lemma}

\begin{proof}
Write $C_N=C_{d+1}(G_N)$, and put $\beta=\pi_{d+1}(\lambda(c))>0$ and
$\rho=\rho(c)$. For fixed $\eta,\gamma>0$, let
$S_{\eta,\gamma}$ be the set of integers $n$ for which there is an
$m$ with $\cG^{(d+1)}_{n,m}\ne\varnothing$ such that
\[
    |2m/n-\rho|<\gamma,
    \qquad
    \PP(U_{n,m}\in\mathcal P)\ge1-\eta,
\]
where $U_{n,m}$ is uniform in $\cG^{(d+1)}_{n,m}$.

We first show that $S_{\eta,\gamma}$ is dense. Otherwise, there
would be a $\delta>0$ and arbitrarily large $A$ with
$S_{\eta,\gamma}\cap(A,(1+\delta)A]=\varnothing$. Choose
\[
    N=\left\lceil\frac{(1+\delta/2)A}{\beta}\right\rceil.
\]
By Theorem~\ref{alt:random-core}, the event
\[
    E=\left\{
      A<|V(C_N)|\le(1+\delta)A,\quad
      \left|\frac{2|E(C_N)|}{|V(C_N)|}-\rho\right|<\gamma
    \right\}
\]
has probability $1-o(1)$ along these values of $A$.
For every pair of vertex and edge counts occurring on $E$, the
conditional probability of $\mathcal P$ is less than $1-\eta$,
by conditional uniformity and the definition of the gap. Therefore
\[
    \PP(C_N\in\mathcal P)
    \le(1-\eta)\PP(E)+\PP(E^c)=1-\eta+o(1),
\]
contrary to the hypothesis.

Apply Lemma~\ref{alt:dense-diagonal} to
$A_t=S_{1/(t+1),1/(t+1)}$. For each term $n_i$ of the resulting
dense sequence, choose an edge count $m_i$ witnessing
$n_i\in A_{t_i}$. Both the density error and the failure probability
are at most $1/(t_i+1)$, and hence tend to zero.
\end{proof}

\begin{lemma}
\label{alt:independent-cores}
There is an $\varepsilon_0>0$ such that, for every
$\varepsilon\in(0,\varepsilon_0)$, there are a dense sequence
$n_j\to\infty$ and edge counts $m_j$ with
$2m_j/n_j\to2d-2\varepsilon$ and
\[
    |\mathcal I^{(d+1)}_{n_j,m_j}|
       =(b_d(2d-2\varepsilon)+o(1))^{n_j}n_j^{m_j}.
\]
\end{lemma}

\begin{proof}
The limiting average degree of the $(d+1)$-core is continuous and
strictly increasing above $\gamma_d$, and equals $2d$ at
$c_d$. Thus, for all sufficiently small $\varepsilon>0$,
there is a $c^*\in(\gamma_d,c_d)$ at which it equals
$2d-2\varepsilon>d+1$.
By Theorem~\ref{alt:independence-threshold},
$G_N\sim G(N,c^*/N)$ and its $(d+1)$-core are a.a.s.\ independent
in the $d$-rigidity matroid. Lemma~\ref{alt:core-transfer} therefore
gives a dense sequence $(n_j,m_j)$ with the stated average degree
and
\[
    \frac{|\mathcal I^{(d+1)}_{n_j,m_j}|}
         {|\cG^{(d+1)}_{n_j,m_j}|}\longrightarrow1.
\]
The claimed count now follows from this ratio and
Theorem~\ref{alt:core-enumeration}.
\end{proof}

We now show that an independent core can be completed to a minimally
$d$-rigid core in many ways. The key input is an upper bound on
the closure of an independent set in the rigidity matroid.

\begin{lemma}
\label{alt:augmentation-count}
Let $n\ge d+1$, let $0\le m<m_n$, and put $r=m_n-m$. Then
\[
    |\cM^{(d+1)}_{n,d}|
    \ge|\mathcal I^{(d+1)}_{n,m}|\,\frac{r!}{(2d^2)^r}.
\]
\end{lemma}

\begin{proof}
Let $\mathcal S$ be the set of sequences $(G_0,\ldots,G_r)$ such that
$G_0\in\mathcal I^{(d+1)}_{n,m}$ and each step adds one edge while
preserving independence in the $d$-rigidity matroid.
Adding edges preserves the minimum degree. Thus $G_r$ is independent
with $m_n$ edges and minimum degree at least $d+1$, that is,
$G_r\in\cM^{(d+1)}_{n,d}$.

We bound $|\mathcal S|$ from both ends. Going forwards, an edge may
be added to $G_t$ precisely when it lies outside
$\operatorname{cl}_d(G_t)$. By Lemma~\ref{alt:closure-bound}, there
are at least
\[
    \frac{r-t}{m_n}\binom n2\ge\frac{n}{2d}(r-t)
\]
choices for this edge. Going backwards, there are at most
$m_n\le dn$ choices for the edge to delete at each step.
Ignoring the minimum-degree constraint during deletion only
increases this upper bound. Hence
\[
    |\mathcal I^{(d+1)}_{n,m}|
       \prod_{t=0}^{r-1}\frac{n(r-t)}{2d}
    \le|\mathcal S|
    \le|\cM^{(d+1)}_{n,d}|(dn)^r.
\]
Dividing gives the claimed bound.
\end{proof}

\begin{lemma}
\label{alt:augmentation-loss}
Suppose $n_j\to\infty$ and $2m_j/n_j\to2d-2\varepsilon$, where
$\varepsilon>0$, and put $r_j=m_{n_j}-m_j$. Then
\[
    |\cM^{(d+1)}_{n_j,d}|
    \ge |\mathcal I^{(d+1)}_{n_j,m_j}|\,n_j^{r_j}
       \exp\!\left[-\bigl(E_d(\varepsilon)+o(1)\bigr)n_j\right],
\]
where $E_d(\varepsilon)=\varepsilon\log(2ed^2/\varepsilon)$ tends
to zero as $\varepsilon\downarrow0$.
\end{lemma}

\begin{proof}
Write $\varepsilon_j=r_j/n_j$, so that
$\varepsilon_j\to\varepsilon$. Stirling's bound gives
\begin{align*}
    \frac{r_j!}{(2d^2)^{r_j}}
    &\ge n_j^{r_j}
       \exp\!\left[n_j\bigl(\varepsilon_j\log\varepsilon_j
                   -\varepsilon_j(1+\log(2d^2))\bigr)\right]\\
    &=n_j^{r_j}
       \exp\!\left[-\bigl(E_d(\varepsilon)+o(1)\bigr)n_j\right].
\end{align*}
Now apply Lemma~\ref{alt:augmentation-count}.
\end{proof}

We are now ready to use diagonalization once more.

\begin{lemma}
\label{alt:optimal-dense-subsequence}
There is a dense sequence $(n_i)$ such that
\[
    a_{n_i}\longrightarrow\log b_d.
\]
\end{lemma}

\begin{proof}
Choose $\varepsilon_t\downarrow0$ with $\varepsilon_t<\varepsilon_0$.
Lemmas~\ref{alt:independent-cores} and~\ref{alt:augmentation-loss}
give, for each $t$, a dense set $A_t$ such that, using
$m_n=dn-O_d(1)$,
\[
    \liminf_{\substack{n\to\infty\\n\in A_t}}a_n
    \ge\log b_d(2d-2\varepsilon_t)-E_d(\varepsilon_t).
\]
Choose $R_t$ so that every $n\in A_t$ with $n>R_t$ satisfies
\[
    a_n\ge
       \log b_d(2d-2\varepsilon_t)-E_d(\varepsilon_t)-1/t.
\]
Lemma~\ref{alt:dense-diagonal} yields a dense sequence
$n_i\in A_{t_i}$ with $n_i>R_{t_i}$ and $t_i\to\infty$.
Continuity of $b_d$ and $E_d(\varepsilon)\to0$ give
$\liminf_i a_{n_i}\ge\log b_d(2d)=\log b_d$.
The upper bound~\eqref{alt:core-upper-bound} gives the matching
upper limit.
\end{proof}

\subsection{Monotonicity}

The subsequence just produced is dense, but it is still a subsequence.
We extend the limit to every $n$ by counting $1$-extensions from
graphs on $n$ vertices to graphs on $n+1$ vertices.
A \emph{$d$-dimensional $1$-extension} deletes an edge $uv$, adds
a new vertex $w$, and joins $w$ to $u$, $v$, and $d-1$ other distinct
vertices. It is known that a $1$-extension preserves minimal
$d$-rigidity~\cite{Henneberg1911,rigidty_general_tibor}. It also
preserves minimum degree at least $d+1$, since no old degree
decreases and the new vertex has degree $d+1$.

\begin{lemma}
\label{alt:entropy-monotonicity}
There is a constant $A_d\ge0$ such that, for all sufficiently large $n$,
\[
    (n+1)a_{n+1}-na_n\ge-A_d.
\]
\end{lemma}

\begin{proof}
Join $G\in\cM^{(d+1)}_{n,d}$ to
$G'\in\cM^{(d+1)}_{n+1,d}$ when $G'$ is obtained from $G$ by a
$1$-extension whose new vertex is labelled $n+1$.
Every resulting graph is minimally $d$-rigid with minimum degree
at least $d+1$.

For each $G$, distinct choices of an edge and a $(d-1)$-subset of
the remaining $n-2$ vertices give distinct graphs $G'$.
Consequently, the number of adjacent pairs is
\[
    m_n\binom{n-2}{d-1}|\cM^{(d+1)}_{n,d}|.
\]
Conversely, a reverse $1$-extension is possible only if the
distinguished vertex $n+1$ has degree $d+1$. There are then at most
$\binom{d+1}{2}$ choices for the edge to be reinserted among its
neighbours. Therefore
\[
    |\cM^{(d+1)}_{n+1,d}|
    \ge
    \frac{m_n\binom{n-2}{d-1}}{\binom{d+1}{2}}
       |\cM^{(d+1)}_{n,d}|.
\]
Lemma~\ref{alt:optimal-dense-subsequence} gives a nonempty family
for at least one sufficiently large $n$. Repeated $1$-extensions
therefore make these families nonempty for every sufficiently large
$n$. Taking logarithms in the last inequality gives
\[
    \log\frac{|\cM^{(d+1)}_{n+1,d}|}{|\cM^{(d+1)}_{n,d}|}
       \ge d\log n-O_d(1).
\]
On the other hand,
\[
    dn\log n-d(n+1)\log(n+1)
    =-d\log n-d(n+1)\log(1+1/n)
    =-d\log n+O_d(1).
\]
By the definition of $a_n$, these two estimates prove the claim.
\end{proof}

The following interpolation argument completes the passage to all sizes.

\begin{lemma}[From a dense subsequence to all sizes]
\label{alt:dense-upgrade}
Let $(a_n)$ be an eventually real sequence satisfying
$(n+1)a_{n+1}-na_n\ge-A$ for some constant $A\ge0$ and all
sufficiently large $n$. If $a_{n_i}\to L$ along a dense sequence,
then $a_n\to L$.
\end{lemma}

\begin{proof}
For $n_i\le n\le n_{i+1}$, summing the increment bound on either
side of $n$ gives
\[
    \frac{n_i}{n}a_{n_i}-A\left(1-\frac{n_i}{n}\right)
    \le a_n\le
    \frac{n_{i+1}}n a_{n_{i+1}}
       +A\left(\frac{n_{i+1}}n-1\right).
\]
Density gives $n_i/n,n_{i+1}/n\to1$, and both bounds tend to $L$.
\end{proof}

\begin{proof}[Proof of Theorem~\ref{thm:rigid-cores}]
Lemmas~\ref{alt:optimal-dense-subsequence} and~\ref{alt:entropy-monotonicity}
verify the hypotheses of Lemma~\ref{alt:dense-upgrade}, so
$a_n\to\log b_d$.
Hence
\[
    \log|\cM^{(d+1)}_{n,d}|
    =dn\log n+n\log b_d+o(n),
\]
which proves the theorem.
Dividing by $|\cG^{(d+1)}_n|$ gives the equivalent probability
statement $(1-o(1))^n$.
\end{proof}

\section{Mantles and parking functions}
\label{sec:attachments}

We turn to the proof of Theorem \ref{thm:acyclic-attachments}.
Parking functions were introduced by Konheim and Weiss~\cite{KonheimWeiss1966},
who showed that there are $n^{n-2}$ parking functions of length $n-1$,
the same as the number of labelled trees on $[n]$ given by Cayley's formula.
A bijective proof was later given by Sch\"utzenberger~\cite{Schutzenberger1968}.
In our terminology, spanning trees are in bijection with $\cD_{k,n-k}^{(0,d)}$ when $d=k=1$, and this bijection extends to larger $d$ and $k$.
For nondecreasing positive integers $a_1,\ldots,a_\ell$, an
\emph{$a$-parking function} is a positive-integer vector $(f_1,...,f_\ell)$ whose
nondecreasing rearrangement $f_{(1)}\le\cdots\le f_{(\ell)}$ satisfies
$f_{(i)}\le a_i$ for every $i$~\cite{KungYan2003}.
Related bijections for acyclic automata appear in
Priez~\cite{Priez2015}.

\begin{lemma}
\label{lem:parking-bijection}
Let $d\le k\le n$ be integers, and set $a_i=\binom{k+i-1}{d}$, for $1 \le i \le n-k$.
Then there is a bijection between $\cD_{k,n-k}^{(0,d)}$ and $a$-parking functions.
\end{lemma}

\begin{proof}

Given $D\in \cD_{k,n-k}^{(0,d)}$, we first construct a canonical ordering $$\pi:[n-k]\to \{k+1,...,n\}\,.$$
Indeed, after $\pi(1),...,\pi(i-1)$ have been chosen, find the minimum index $v\in \{k+1,...,n\}$ that has not been chosen yet and whose $d$ out-neighbors $N_v$ all belong to $[k]\cup\pi([i-1])$. Such a vertex $v$ exists by the assumption that $D$ is acyclic. Set $\pi(i):=v$. In addition, consider colex order on $d$-subsets of $[k]\cup\pi([i-1])$ induced by the order $1,..,k,\pi(1),...,\pi(i-1)$, and let $1\le f_{v-k} \le a_i$ be the index of $N_v$ in this colex order. By construction, the resulting vector $(f_1,...,f_{n-k})$ is an $a$-parking function.

On the other hand, given an $a$-parking function $(f_1,...,f_{n-k})$, we can construct a unique $D\in \cD_{k,n-k}^{(0,d)}$ as follows. We first construct the canonical ordering $\pi$. Indeed, after $\pi(1),...,\pi(i-1)$ have been chosen, find the minimum index $v\in \{k+1,...,n\}$ that has not been chosen yet such that $f_{v-k}\le a_i$, and set the out-neighbors of $v$ to be the $f_{v-k}$-th $d$-subset in colex order of $[k]\cup\pi([i-1])$. Such a $v$ exists by the assumption that the given sequence is an $a$-parking function. By construction, the resulting digraph is acyclic and has the required out-degrees.

Colex order lists the $d$-subsets of any initial segment first and preserves
their indices when the ground set is extended.

To verify that these maps are inverse, let $\pi_D$ and $\pi_f$ denote the canonical orderings obtained from $D$ and $f$, respectively. We claim that $\pi_D=\pi_f$ whenever $D$ and $f$ are related by either construction. Indeed, suppose inductively that the two orderings agree up to step $i-1$. By the definition of colex order, a remaining vertex $v$ satisfies $N_v\subseteq [k]\cup\pi([i-1])$ if and only if $f_{v-k}\le a_i$. Hence, both constructions choose the same vertex at step $i$, proving the claim. Finally, once $\pi$ is fixed, the out-neighbor sets $(N_v)$ and their colex indices $(f_{v-k})$ uniquely determine each other. This proves the bijection.
\end{proof}

\begin{figure}[tbp]
    \centering
    \includegraphics[width=\textwidth]{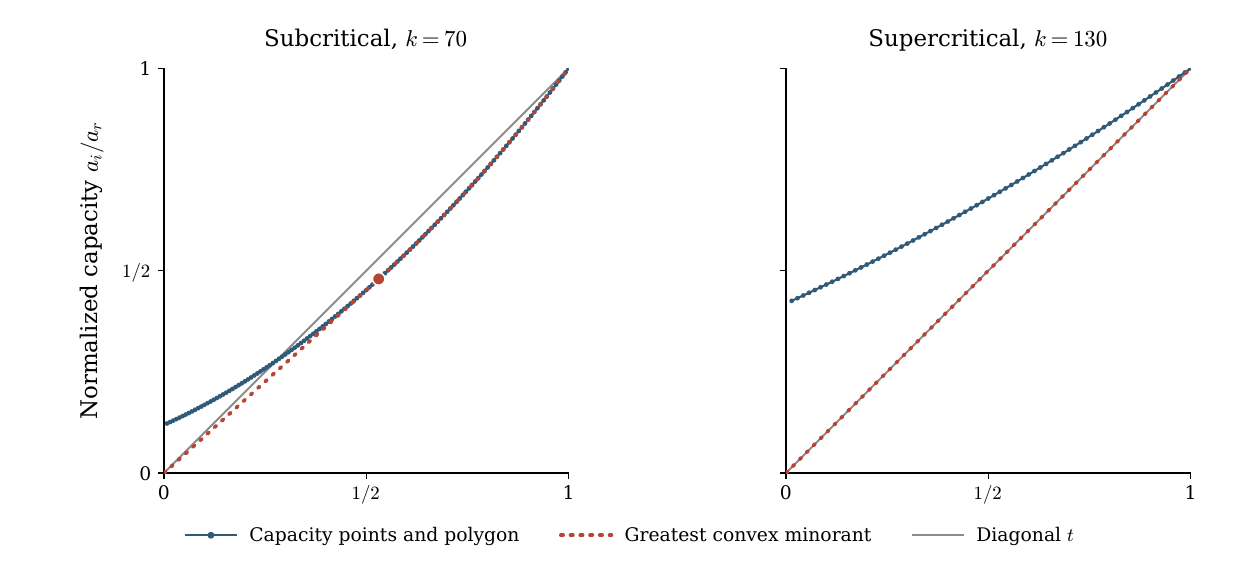}
    \caption{Normalized capacity points $(i/r,a_i/a_r)$ and their convex
    minorants for $d=2$, $n=200$, and $k=70$ (left) and $k=130$ (right).
    These values lie below and above the critical fraction $k/n=1/2$.
    In the left panel, the highlighted point marks the end of the initial
    chord, after which the minorant follows the capacity polygon;
    in the right panel, it is the diagonal.}
    \label{fig:parking-minorants}
\end{figure}

We now asymptotically count $a$-parking functions with
\[
    a_i=\binom{k+i-1}{d},\qquad 1\le i\le r:=n-k.
\]
Suppose that $k/n\to\alpha\in(0,1)$ and define
\[
    \rho(t)=\bigl(\alpha+(1-\alpha)t\bigr)^d,
    \qquad 0\le t\le1.
\]
Since $a_r=\binom{n-1}{d}$, the binomial formula gives
$a_i/a_r=\rho(i/r)+o(1)$, uniformly for $1\le i\le r$.
Let $g$ be the greatest convex minorant of this curve together with
the origin. It is the limit of the convex minorants of the polygons through
\[
    (0,0),\qquad
    \left(\frac{i}{r},\frac{a_i}{a_r}\right),
    \quad 1\le i\le r.
\]
The origin is included even though $\rho(0)=\alpha^d>0$.
If $\alpha\ge1-1/d$, then $\rho(t)\ge t$ and $g(t)=t$.
Otherwise, the line from the origin is tangent to $\rho$ at
\[
    s=\frac{\alpha}{(d-1)(1-\alpha)},
    \qquad \rho(s)=s\rho'(s),
\]
and
\[
    g(t)=
    \begin{cases}
        \rho'(s)t,&0\le t\le s,\\
        \rho(t),&s\le t\le1.
    \end{cases}
\]
Figure~\ref{fig:parking-minorants}
illustrates the subcritical and supercritical regimes for $d=2$.

\begin{lemma}
\label{lem:parking-convex-minorant}
Fix $d\ge2$ and suppose that $k/n\to\alpha\in(0,1)$. With
$r$, $a_i$, and $g$ as above, let $f_1,\ldots,f_r$ be i.i.d.\ uniform
on $\left[\binom{n-1}{d}\right]$. Then
\begin{equation}
    \PP\bigl((f_i)_{i=1}^{r}\text{ is an }a\text{-parking function}\bigr)
    =\exp\left(r\int_0^1\log g'(t)\,dt+o(n)\right).
    \label{eq:parking-probability-rate}
\end{equation}
\end{lemma}

\begin{proof}
Let $E$ be the parking event, and let $F=(F_1,\ldots,F_r)$ have
the conditional law of $(f_1,\ldots,f_r)$ given $E$. Write $H$ for
Shannon entropy with natural logarithms, and recall that
$a_r=\binom{n-1}{d}$.
Since the coordinates of $F$ are exchangeable, the parking condition gives
\[
    \PP(F_1\le a_i)
    =\frac1r\mathbb{E}\bigl|\{j:F_j\le a_i\}\bigr|
    \ge\frac{i}{r},
    \qquad 1\le i\le r.
\]
Let $\nu$ maximize entropy among probability distributions on $[a_r]$
satisfying these cumulative constraints, and write $q(x)=\nu(\{x\})$.
Since $F$ is uniform among the $a$-parking functions, entropy
subadditivity yields
\begin{equation}
\label{eq:parking-entropy-bound}
    \log\PP(E)
    =H(F)-r\log a_r
    \le r\bigl(H(\nu)-\log a_r\bigr).
\end{equation}

Set $a_0=0$. Averaging the probabilities within each interval
$(a_{i-1},a_i]$ preserves all the constraints and increases entropy
unless they are already equal. Thus $q(x)=q_i$ throughout this interval.
Moving mass from a later atom to an earlier one shows that
$q_1\ge\cdots\ge q_r$: otherwise such a transfer would preserve the
constraints and increase entropy.
Moreover, if $q_i>q_{i+1}$, the constraint at $a_i$ must be tight;
otherwise a small transfer of mass from $a_i$ to $a_i+1$ would preserve
the constraints and increase entropy.

Let
\[
    0=m_0<m_1<\cdots<m_L=r
\]
be the endpoints of the maximal constant blocks of $(q_i)_{i=1}^r$.
The cumulative constraints are equalities at these endpoints, so
\[
    (a_{m_j}-a_{m_{j-1}})q_{m_j}
    =\frac{m_j-m_{j-1}}{r},
    \qquad 1\le j\le L.
\]
Let $g_n$ be the polygonal line joining
\[
    \left(\frac{m_j}{r},\frac{a_{m_j}}{a_r}\right),
    \qquad 0\le j\le L.
\]
Its slope on the $j$-th block is $1/(a_rq_{m_j})$, so $g_n$ is convex.
For $m_{j-1}<i<m_j$, feasibility and tightness at the left endpoint give
\[
    \frac{i-m_{j-1}}{r}
    \le (a_i-a_{m_{j-1}})q_{m_j},
\]
which is precisely $g_n(i/r)\le a_i/a_r$.
Thus $g_n$ lies below all the prescribed points.
Every convex minorant lies below the chord joining each pair of
consecutive block endpoints, and hence below $g_n$.
Therefore $g_n$ is the greatest convex minorant of $(0,0)$ and
$(i/r,a_i/a_r)$, $1\le i\le r$.
The block masses also give
\begin{equation}
\label{eq:parking-entropy-integral}
    H(\nu)-\log a_r
    =\sum_{j=1}^L\frac{m_j-m_{j-1}}{r}
        \log\frac{1}{a_rq_{m_j}}
    =\int_0^1\log g_n'(t)\,dt.
\end{equation}

We next establish bounds that will be used both for concentration
and for passage to the limit.
Since $a_1q_1\le1$, we have $q(x)\le1/a_1$.
For the lower bound, if $\nu$ is uniform then $q(x)=1/a_r$.
Otherwise, writing $m=m_{L-1}\ge1$, we have
\[
    q_r=\frac{r-m}{r(a_r-a_m)}
    \ge\frac{1}{r\binom{n-2}{d-1}}
    =\frac{n-1}{dr\,a_r}
    \ge\frac{1}{d\,a_r},
\]
where we used
$a_i-a_{i-1}\le\binom{n-2}{d-1}$ for $i>m$ and $r\le n-1$.
Consequently,
\[
    \frac{1}{d\,a_r}\le q(x)\le\frac{1}{a_1}.
\]
In particular,
\begin{equation}
\label{eq:parking-uniform-bounds}
    \frac{a_1}{a_r}\le g_n'(t)\le d,
    \qquad
    \log\frac{a_1}{a_r}
    \le-\log q(x)-\log a_r\le\log d.
\end{equation}
Since $a_1/a_r\to\alpha^d>0$, these bounds are uniform for all
sufficiently large $n$.

For the lower bound on $\PP(E)$, let $\nu_0$ be uniform on
$\{a_1,\ldots,a_r\}$.
The cumulative constraints imply that $\nu$ is stochastically
dominated by $\nu_0$: writing $F_\nu$ and $F_{\nu_0}$ for their
distribution functions, we have
\[
    F_\nu(x)\ge F_\nu(a_i)\ge\frac{i}{r}=F_{\nu_0}(x),
    \qquad a_i\le x<a_{i+1},\quad 0\le i<r.
\]
Since decreasing any coordinate preserves the parking condition,
\[
    \PP_{\nu^{\otimes r}}(E)
    \ge\PP_{\nu_0^{\otimes r}}(E)
    =\frac{(r+1)^{r-1}}{r^r}
    \ge\frac1r.
\]
Here the equality follows from the classical parking-function
count~\cite{KonheimWeiss1966}, after replacing each $a_i$ by $i$.

Let $Z_1,\ldots,Z_r$ be independent with law $\nu$, and set
\[
    S=\sum_{j=1}^r-\log q(Z_j).
\]
Its mean is $rH(\nu)$.
By \eqref{eq:parking-uniform-bounds}, Hoeffding's inequality gives,
for some $C=C(\alpha,d)$ and all sufficiently large $n$,
\[
    \PP_{\nu^{\otimes r}}\left(
        S<rH(\nu)-C\sqrt{r\log r}
    \right)\le r^{-2}.
\]
Thus the parking vectors satisfying
$S\ge rH(\nu)-C\sqrt{r\log r}$ have total probability at least $1/(2r)$.
Each such vector has probability at most
$\exp(-rH(\nu)+C\sqrt{r\log r})$, so there are at least
\[
    \frac{1}{2r}
    \exp\left(rH(\nu)-C\sqrt{r\log r}\right)
\]
$a$-parking functions.
Together with \eqref{eq:parking-entropy-bound}, this yields
\begin{equation}
\label{eq:parking-discrete-rate}
    \log\PP(E)
    =r\bigl(H(\nu)-\log a_r\bigr)
      +O_{\alpha,d}\left(\sqrt{r\log r}\right).
\end{equation}

Finally, we pass to the limit.
Set
\[
    b(t)=(\alpha+(1-\alpha)t)^d,
    \qquad
    \delta_n=\max_{1\le i\le r}
    \left|\frac{a_i}{a_r}-b(i/r)\right|.
\]
The binomial formula for $a_i$ gives $\delta_n\to0$.
Recall that $g$ is the greatest convex minorant of $b$ together
with the point $(0,0)$.
Since $g-\delta_n$ lies below all the discrete points, maximality gives
$g-\delta_n\le g_n$.
Conversely, monotonicity of $g_n$ and the fact that $b$ is
$d$-Lipschitz imply
\[
    g_n(t)\le b(t)+\delta_n+\frac dr,
    \qquad 0\le t\le1,
\]
by comparison with the next mesh point.
Hence $g_n-\delta_n-d/r$ is a convex minorant of $b$ and $(0,0)$,
so
\[
    g-\delta_n\le g_n\le g+\delta_n+\frac dr.
\]
Thus $g_n\to g$ uniformly.
Convexity now implies $g_n'\to g'$ almost everywhere, and
\eqref{eq:parking-uniform-bounds} permits dominated convergence.
Using \eqref{eq:parking-entropy-integral}, we obtain
\[
    H(\nu)-\log a_r
    =\int_0^1\log g_n'(t)\,dt
    \longrightarrow
    \int_0^1\log g'(t)\,dt.
\]
Since $r=n-k=\Theta(n)$, \eqref{eq:parking-discrete-rate} gives
\[
    \log\PP(E)
    =r\int_0^1\log g'(t)\,dt+o(n),
\]
as required.
\end{proof}

\begin{proof}[Proof of Theorem~\ref{thm:acyclic-attachments}]
First suppose $k/n\to\alpha\in(0,1)$ and put $r=n-k$.
By Lemmas~\ref{lem:parking-bijection} and~\ref{lem:parking-convex-minorant},
\[
    |\cD_{k,n-k}^{(0,d)}|=\binom{n-1}{d}^{r}
    \exp\left(r\int_0^1\log g'(t)\,dt+o(n)\right).
\]
Since
\[
    \frac{\binom{n-1}{d}^{r}}{\binom nd^{r}}
    =\left(1-\frac dn\right)^r=\exp(O_d(1)),
\]
and $r/n\to1-\alpha$, the cost in the theorem is
\begin{equation}
    J_d(\alpha)
    =-(1-\alpha)\int_0^1\log g'(t)\,dt.
    \label{eq:acyclicity-convex-minorant}
\end{equation}
It is zero when $g(t)=t$. Below the threshold, substituting the
two parts of $g$ gives
\[
    J_d(\alpha)=-(1-\alpha)
    \left(s\log\rho'(s)+\int_s^1\log\rho'(t)\,dt\right).
\]
Elementary integration yields the stated closed form.
To obtain the uniform error-term, write
\[
    p_n(k)=\frac{|\cD_{k,n-k}^{(0,d)}|}{\binom nd^{n-k}}.
\]
This is nondecreasing in $k$. Indeed, deleting the outgoing edges
of $k'-k$ vertices enlarges the sink set from size $k$ to size $k'$.
Each resulting digraph has at most $\binom{n-1}{d}^{k'-k}$ preimages,
and $\binom{n-1}{d}\le\binom nd$.

For $k=d$, adding a clique on the sink set and forgetting the
orientations gives a bijection with the graphs grown from that fixed
seed by $0$-extensions, by the uniqueness of the acyclic orientation.
Allowing all seed sets changes the count
by at most a factor of $\binom nd$. Thus \eqref{eq:zero-extension-count}
gives
\[
    \frac{1}{n}\log p_n(d)
    =-(d-1)-\log((d-1)!)+\log(d!)+o(1)
    =-J_d(0)+o(1).
\]
At the other endpoint, $p_n(n)=1$.

Fix a mesh of $[0,1]$ on which $J_d$ varies by at most $\varepsilon$
on each interval. At each interior mesh point $t$, apply the estimate
above with $k=\lfloor tn\rfloor$. At the endpoints use $k=d$ and $k=n$.
Monotonicity bounds $n^{-1}\log p_n(k)$ between its values at the
adjacent mesh indices. By continuity of $J_d$, these bounds differ
from $-J_d(k/n)$ by at most $\varepsilon+o(1)$ uniformly in $k$.
Letting $\varepsilon\downarrow0$ proves the theorem.
\end{proof}

\section{Proofs of the main results}
\label{sec:main-proofs}
We start with a proof of \eqref{eq:rigid-core-decomposition}. 
Every minimally $d$-rigid graph with more than $d$ vertices has minimum
degree at least $d$. Removing a vertex of degree $d$ leaves a minimally $d$-rigid graph. Repeatedly remove such vertices while
more than $d$ vertices remain. We end either with $K_d$, in which case
$G\in\cM_{n,d}(0)$, or with a graph of minimum degree at least $d+1$.
In the latter case, the remaining graph is exactly $C_{d+1}(G)$ and is
minimally $d$-rigid. Write $C$ for its vertex set and $k=|C|$; necessarily
$k\ge d+2$.

Order the removed vertices by the time they were peeled. Direct each
mantle edge from the vertex removed earlier to the vertex removed later,
or to its endpoint in $C$. Every removed vertex has exactly $d$
outgoing edges, every core vertex has out-degree $0$, and the orientation
is acyclic. After relabeling $C$ as $[k]$, this gives a member of
$\cD_{k,n-k}^{(0,d)}$. The orientation is unique: fix a peeling order.
Inductively, edges to earlier vertices already point into the current
vertex, so all its remaining $d$ edges must point outward.

Conversely, choose a minimally $d$-rigid core on $C$ and any digraph
in $\cD_{k,n-k}^{(0,d)}$, with sink set $C$. Add the other vertices in reverse
topological order, joining each to its $d$ out-neighbors already present.
These are $0$-extensions, so the resulting graph is minimally $d$-rigid.
Peeling the added vertices in topological order recovers precisely the
chosen core. Together with uniqueness of the orientation, this gives a
bijection between graphs with nonempty core and the choices of the core
vertex set, the core graph, and the mantle digraph. The number with core
size $k$ is therefore
$\binom nk|\cM_{k,d}^{(d+1)}|\,|\cD_{k,n-k}^{(0,d)}|$.
Summing over $k$ and adding the empty-core family $\cM_{n,d}(0)$ proves
\eqref{eq:rigid-core-decomposition}.

\begin{proof}[Proof of Theorem~\ref{thm:enumeration}]

Write $\cM_{n,d}(k)$ for the graphs in $\cM_{n,d}$ whose $(d+1)$-core
has $k$ vertices. For $k>0$ the corresponding summand in
\eqref{eq:rigid-core-decomposition} gives
\[
    |\cM_{n,d}(k)|=\binom nk|\cM_{k,d}^{(d+1)}|\,|\cD_{k,n-k}^{(0,d)}|.
\]
Put $\alpha=k/n$. By Theorems~\ref{thm:rigid-cores}
and~\ref{thm:acyclic-attachments} and Stirling's formula, for every
feasible $k>0$ we have
\begin{equation}
    \label{eq:rigid-slice-rate}
    \begin{aligned}
    &\frac1n\log|\cM_{n,d}(k)|-d\log n
    =H(\alpha)+\alpha\log(b_d\alpha^d)
      -(1-\alpha)\log(d!)-J_d(\alpha)+o(1)\\
    &\qquad=o(1)+\begin{cases}
        -(d-1)-\log((d-1)!)
        +\alpha\log\!\left(\dfrac{e^d b_d\,d!(d-1)^{d-1}}{d^d}\right),
        &0\le\alpha\le1-\dfrac1d,\\[8pt]
        \alpha\log b_d+(d-1)\alpha\log\alpha
        -(1-\alpha)\log\bigl(d!(1-\alpha)\bigr),
        &1-\dfrac1d\le\alpha\le1,
    \end{cases}
    .
    \end{aligned}
\end{equation}
Here $H(\alpha)=-\alpha\log\alpha-(1-\alpha)\log(1-\alpha)$
is the binary entropy, and the second equality follows by direct
calculation. The $o(1)$ error is uniform in $k$, meaning that the
largest absolute error over all feasible $k$ tends to zero as
$n\to\infty$. The core enumeration has logarithmic error $o(k)$ for
large $k$, while the finitely many smaller feasible sizes have bounded
errors. Dividing by $n$ gives $o(1)$ uniformly for $k\le n$.
Write $\widetilde F_d(\alpha)$ for the piecewise expression
above, with $0\log0=0$. At $k=0$, \eqref{eq:zero-extension-count}
gives the same estimate, since
$\widetilde F_d(0)=-(d-1)-\log((d-1)!)$.
Figure~\ref{fig:core-size-rates} shows these rates for $d=2,3$.

\begin{figure}[tbp]
    \centering
    \includegraphics[width=\textwidth]{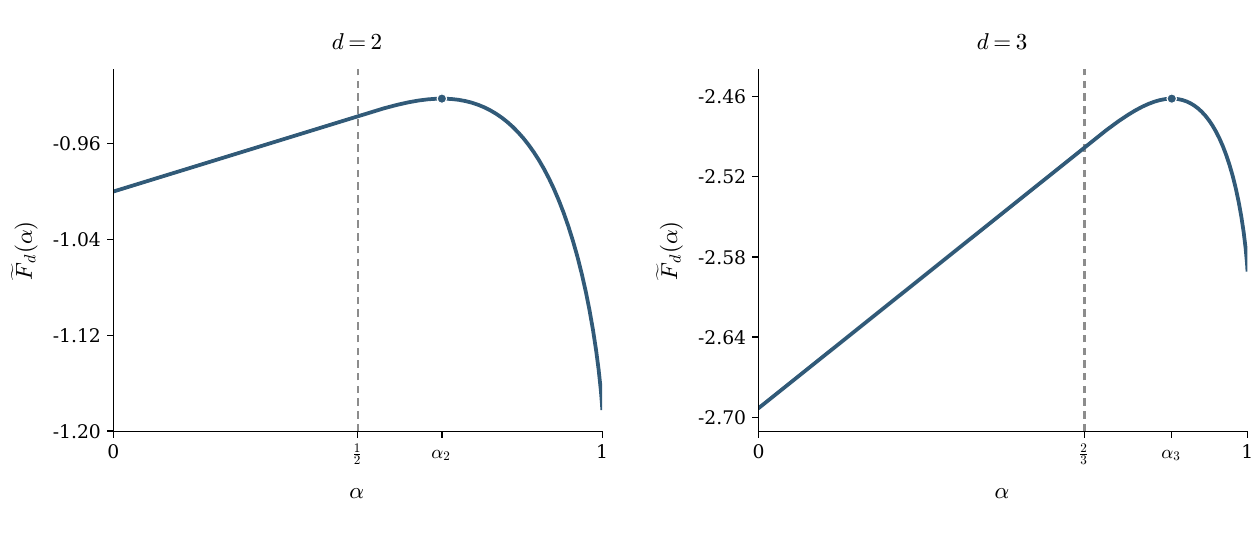}
    \caption{The rate functions $\widetilde F_2$ and $\widetilde F_3$.
    The dashed lines mark $1-1/d$, below which the functions are affine.
    The dots mark their unique maxima at $\alpha_d$.}
    \label{fig:core-size-rates}
\end{figure}

It will be useful to express $b_d$ in terms of $\alpha_d$. By \eqref{eq:poisson-parameter},
$\lambda_d^*=2d\alpha_d$, and therefore
\begin{equation}
    b_d=\frac{e^{-d}}{d!\,\alpha_d^{d-1}(1-\alpha_d)}.
    \label{eq:bd-alpha}
\end{equation}

In addition, we note that $1-1/d<\alpha_d<1$ for all $d\ge2$ since,
for $Z\sim\Poi(\lambda_d^*)$, we have
$\alpha_d=\PP(Z\ge d+1\mid Z\ge d)<1$, and
\[
    0<\operatorname{Var}(Z\mid Z\ge d+1)
    =d(\lambda_d^*-2d+2)
    =2d^2\left(\alpha_d-1+\frac1d\right).
\]
The variance identity follows from the Poisson recurrence and
\eqref{eq:poisson-parameter}.

Observe that $\widetilde F_d$ is linear on $[0,1-1/d]$ and its
slope is positive since, using~\eqref{eq:bd-alpha}, we have
\[
    \log\frac{e^d b_d\,d!(d-1)^{d-1}}{d^d}
    =\log\frac{(1-1/d)^{d-1}/d}
                  {\alpha_d^{d-1}(1-\alpha_d)}>0,
\]
as $\alpha_d>1-1/d$ and $x^{d-1}(1-x)$ is strictly decreasing on
$[1-1/d,1]$.

In addition, $\widetilde F_d$ is strictly concave on $[1-1/d,1]$ since
\[
    \widetilde F_d''(\alpha)
    =\frac{d-1}{\alpha}-\frac1{1-\alpha}
    =\frac{d-1-d\alpha}{\alpha(1-\alpha)}<0,
    \qquad 1-\frac1d<\alpha<1.
\]

Finally, using~\eqref{eq:bd-alpha},
\[
    \widetilde F_d'(\alpha_d)
    =\log\!\left(e^d b_d\,d!\,\alpha_d^{d-1}(1-\alpha_d)\right)=0.
\]

Thus $\alpha_d$ is the unique maximizer on $[0,1]$. Using~\eqref{eq:bd-alpha} gives
\begin{equation}
    \max_{\alpha\in[0,1]}\widetilde F_d(\alpha)
    =\widetilde F_d(\alpha_d)
    =-d\alpha_d-\log(1-\alpha_d)-\log(d!)
    =\log\kappa_d.
    \label{eq:core-profile-maximum}
\end{equation}

Summing the counts in \eqref{eq:rigid-slice-rate} over the at most $n+1$ feasible
core sizes and using \eqref{eq:core-profile-maximum} gives
$|\cM_{n,d}|\le n^{dn}\exp(n\log\kappa_d+o(n))$.
The matching lower bound follows by taking $k=\lfloor\alpha_d n\rfloor$
in \eqref{eq:rigid-slice-rate}.
\end{proof}

\begin{proof}[Proof of Corollary~\ref{cor:rigid-core-concentration}]

Fix $\varepsilon>0$. By continuity and the unique maximum
\eqref{eq:core-profile-maximum}, there is $\gamma>0$ such that
\[
    \widetilde F_d(\alpha)\le\log\kappa_d-2\gamma
    \quad\text{whenever }|\alpha-\alpha_d|\ge\varepsilon.
\]
Summing the counts in \eqref{eq:rigid-slice-rate} over these core sizes
and dividing by $|\cM_{n,d}|$ gives
\[
    \PP\!\left(
        \left|\,|V(C_{d+1}(G))|-\alpha_d n\,\right|>\varepsilon n
    \right)
    \le (n+1)\exp(-2\gamma n+o(n))
    \le e^{-\gamma n}
\]
for all sufficiently large $n$.
\end{proof}

We also find additional structure in a typical minimally $d$-rigid graph
by studying the levels of its mantle.

Let $G\in\cM_{n,d}$. Starting from the core vertex set
$S_0=V(C_{d+1}(G))$, define
\[
    S_{i+1}=S_i\cup\{v\notin S_i:|N_G(v)\cap S_i|\ge d\}.
\]
With $S_{-1}=\varnothing$, the $i$th level is $L_i=S_i\setminus S_{i-1}$.

\begin{theorem}\label{thm:level-concentration}
Set $\beta_0=\sigma_0=\alpha_d$ and $\sigma_{-1}=0$. For $i\ge1$, define
\begin{equation}
    \beta_i=(1-\alpha_d)(\sigma_{i-1}^d-\sigma_{i-2}^d),
    \qquad \sigma_i=\sigma_{i-1}+\beta_i.
    \label{eq:level-profile}
\end{equation}
If $G$ is uniform in $\cM_{n,d}$, then for every $\varepsilon>0$
there is $c=c(d,\varepsilon)>0$ such that, for all sufficiently large $n$,
\[
    \PP\!\left(\sup_{i\ge0}
        \left|\frac{|L_i|}{n}-\beta_i\right|>\varepsilon\right)
    \le e^{-cn}.
\]
\end{theorem}

The limiting level densities sum to one, that is, $\sigma_i\to1$.
Indeed, by \eqref{eq:level-profile},
\[\sigma_i=\sum_{j=0}^i\beta_j
=\alpha_d+(1-\alpha_d)\sum_{j=1}^i(\sigma_{j-1}^d-\sigma_{j-2}^d)
=\alpha_d+(1-\alpha_d)\sigma_{i-1}^d\,,\]
Therefore, $0\le\sigma_i\le1$
by induction, thus $1-\sigma_{i-1}^d\le d(1-\sigma_{i-1})$ and
\[
    1-\sigma_i=(1-\alpha_d)(1-\sigma_{i-1}^d)
        \le d(1-\alpha_d)(1-\sigma_{i-1}).
\]
Since $\alpha_d>1-1/d$, $1-\sigma_i$ tends to zero exponentially.

\begin{table}[htbp]
    \centering
    \begin{tabular}{crrrrrrr}
        \toprule
        $d$ & $\beta_0$ & $\beta_1$ & $\beta_2$ & $\beta_3$
            & $\beta_4$ & $\beta_5$ & $\beta_6$ \\
        \midrule
        $2$ & 0.672000 & 0.148120 & 0.072492 & 0.040724
            & 0.024390 & 0.015128 & 0.009580 \\
        $3$ & 0.845245 & 0.093453 & 0.034551 & 0.014661
            & 0.006545 & 0.002985 & 0.001375 \\
        \bottomrule
    \end{tabular}
    \caption{Limiting level densities for $d=2,3$, rounded to six decimal places.}
    \label{tab:level-densities}
\end{table}

\begin{proof}
Fix a level profile $\ell=(\ell_0,\ldots,\ell_m)$ of positive integers
summing to $n$, and set $\ell_i=0$ for $i>m$. For any fixed $\eta>0$,
we may assume $|\ell_0/n-\beta_0|\le\eta$ by
Corollary~\ref{cor:rigid-core-concentration}: the probability of all other
core sizes is exponentially small. Put $s_{-1}=0$ and
$s_i=\ell_0+\cdots+\ell_i$. The number of graphs in $\cM_{n,d}$ with profile $\ell$ is
\[
   N_{\ell}:= \binom{n}{\ell_0,\ldots,\ell_m}
    |\cM_{\ell_0,d}^{(d+1)}|
    \prod_{i=1}^m
    \left[\binom{s_{i-1}}d-\binom{s_{i-2}}d\right]^{\ell_i}.
\]
Indeed, choose the level sets and a minimally rigid core on $L_0$.
Inductively, no mantle edge is directed from $S_{i-1}$ to its complement,
so every edge joining a remaining vertex to $S_{i-1}$ points into $S_{i-1}$.
A vertex belongs to $L_i$ precisely when its $d$ out-neighbours lie in
$S_{i-1}$ and at least one lies in $L_{i-1}$. Every such choice gives a
graph with profile $\ell$, and its unique mantle orientation recovers the choices.

Let $T\ge1$ be chosen later. We now count all graphs whose level sizes
agree with $\ell$ through level $T$, regardless of the later levels.
Put $r=n-s_T$ and let
\[
    p_i=\frac{\binom{s_{i-1}}d-\binom{s_{i-2}}d}{\binom nd}
    \quad(1\le i\le T),\qquad
    p_{T+1}=1-\frac{\binom{s_{T-1}}d}{\binom nd}.
\]
Each of the $r$ remaining vertices must choose a $d$-set not contained
in $S_{T-1}$, since otherwise it would belong to $S_T$. Thus
\begin{align*}
    \sum_{\substack{\ell':\,\ell'_i=\ell_i\\0\le i\le T}} N_{\ell'}
    &\le \binom{n}{\ell_0,\ldots,\ell_T,r}|\cM_{\ell_0,d}^{(d+1)}|
         \binom nd^{n-\ell_0}\left(\prod_{i=1}^T p_i^{\ell_i}\right)p_{T+1}^r\\
    &= \left[\binom n{\ell_0}|\cM_{\ell_0,d}^{(d+1)}|\binom nd^{n-\ell_0}\right]
       \cdot\binom{n-\ell_0}{\ell_1,\ldots,\ell_T,r}
       \left(\prod_{i=1}^T p_i^{\ell_i}\right)p_{T+1}^r\\
    &\le n^{dn}(\kappa_d+o(1))^n
       \cdot\binom{n-\ell_0}{\ell_1,\ldots,\ell_T,r}
       \left(\prod_{i=1}^T p_i^{\ell_i}\right)p_{T+1}^r .
\end{align*}
Indeed, $J_d(\ell_0/n)=0$ for $\ell_0/n$ sufficiently close to $\beta_0>1-1/d$.
By Theorem~\ref{thm:acyclic-attachments}, $\binom nd^{n-\ell_0}$ then counts
the mantles up to a factor of $\exp(o(n))$. Thus the bracketed factor is
$e^{o(n)}|\cM_{n,d}(\ell_0)|$ uniformly over these core sizes. Since
$\cM_{n,d}(\ell_0)\subseteq\cM_{n,d}$, Theorem~\ref{thm:enumeration} gives the last bound.

Let $(X_1,\ldots,X_{T+1})$ have the multinomial distribution with
$n-\ell_0$ trials and outcome probabilities $p_1,\ldots,p_{T+1}$.
Note that
\[
    \binom{n-\ell_0}{\ell_1,\ldots,\ell_T,r}
    \left(\prod_{i=1}^T p_i^{\ell_i}\right)p_{T+1}^r
    =\PP\bigl((X_1,\ldots,X_{T+1})=(\ell_1,\ldots,\ell_T,r)\bigr).
\]
Fix $\delta>0$, and suppose that $|\ell_i/(n-\ell_0)-p_i|\ge\delta$ for some $1\le i\le T$.
Since $X_i$ has the binomial distribution with parameters $n-\ell_0$ and $p_i$,
Hoeffding's inequality gives
\[
     \binom{n-\ell_0}{\ell_1,\ldots,\ell_T,r}
    \left(\prod_{i=1}^T p_i^{\ell_i}\right)p_{T+1}^r
    \le \PP\!\left(\left|\frac{X_i}{n-\ell_0}-p_i\right|\ge\delta\right)
    \le 2e^{-2\delta^2(n-\ell_0)}.
\]

There are at most $(n+1)^{T+1}$ choices of $\ell_0,\ldots,\ell_T$, and
$n-\ell_0\ge(1-\beta_0-\eta)n$. Summing the bounds above, the fraction of
graphs in $\cM_{n,d}$ with such a deviation and $|\ell_0/n-\beta_0|\le\eta$
is at most $\exp(-2\delta^2(1-\beta_0-\eta)n+o(n))$.

\medskip

It therefore remains to prove the following deterministic claim.
For every $\varepsilon>0$ there exist an integer $T\ge1$ and
$\eta,\delta>0$ such that, for all sufficiently large $n$, every level
profile satisfying
\[
    \left|\frac{\ell_0}{n}-\beta_0\right|\le\eta
    \qquad\text{and}\qquad
    \left|\frac{\ell_i}{n-\ell_0}-p_i\right|<\delta
    \quad(1\le i\le T)
\]
    also satisfies $|\ell_i/n-\beta_i|<\varepsilon$ for \emph{every} $i\ge0$.

To prove the claim, choose $T$ so that $1-\sigma_T<\varepsilon/4$.
Summing the first $i$ discrepancies and using the definition of $p_j$ gives
\[
    \left|\frac{s_i-\ell_0}{n-\ell_0}
          -\frac{\binom{s_{i-1}}d}{\binom nd}\right|
    =\left|\sum_{j=1}^i\left(\frac{\ell_j}{n-\ell_0}-p_j\right)\right|
    \le i\delta,\qquad 1\le i\le T.
\]
We next compare $s_i/n$ with $\sigma_i$ through this fixed level $T$.
Multiplying the bound on the first $i$ discrepancies by $(n-\ell_0)/n$
and using $\binom{s_{i-1}}d/\binom nd=(s_{i-1}/n)^d+O_d(1/n)$ gives
\[
    \frac{s_i}{n}=\frac{\ell_0}{n}
       +\left(1-\frac{\ell_0}{n}\right)\left(\frac{s_{i-1}}n\right)^d
       +O_d(T\delta+1/n).
\]
Replacing $\ell_0/n$ by $\alpha_d$ changes the right-hand side
by at most $\eta$, since the coefficient of $\ell_0/n$ is
$1-(s_{i-1}/n)^d\in[0,1]$. An error in the previous cumulative size
does not grow under the recurrence: for $x,y\in[0,1]$,
\[
    (1-\alpha_d)|x^d-y^d|\le d(1-\alpha_d)|x-y|\le |x-y|.
\]
The initial error $|s_0/n-\sigma_0|$ is at most $\eta$, and each of the
next $T$ steps adds at most $\eta+O_d(T\delta+1/n)$ to the error.
Since $T$ is fixed, we can choose $\eta,\delta$ small enough, and then
$n$ large enough, that $|s_i/n-\sigma_i|<\varepsilon/4$ for every $i\le T$.

For $1\le i\le T$, the size of level $i$ is the difference of two
consecutive cumulative sizes. Hence
\[
    \left|\frac{\ell_i}{n}-\beta_i\right|
    \le \left|\frac{s_i}{n}-\sigma_i\right|
       +\left|\frac{s_{i-1}}n-\sigma_{i-1}\right|<\varepsilon/2.
\]
Level zero is already covered since $s_0=\ell_0$. For the later levels
it suffices to bound their total mass. The limiting profile has less
than $\varepsilon/4$ left after $T$, and the observed profile has at most
$1-s_T/n\le 1-\sigma_T+|s_T/n-\sigma_T|<\varepsilon/2$ left.
Thus for every $i>T$, both $\ell_i/n$ and $\beta_i$ lie in
$[0,\varepsilon/2)$, which proves the claim.
\end{proof}

\subsubsection*{Declaration of AI use}

The main substantive contribution of AI (ChatGPT Sol 5.6) to this work
was suggesting the bijection between mantles and generalized parking
functions, which became a key ingredient of our proof. Before this
suggestion, we had obtained a different proof of
Theorem~\ref{thm:enumeration}, based on the fact that the
$((d+1)+d)$-core of $G(n,p)$, conditional on its size and number of
edges, is uniformly distributed. 
This proof appears in the Master's thesis of the first author.
However, this approach did not yield
exponential concentration of the core size. 
Beyond this mathematical
contribution, we used AI to improve the writing, clarity, and
presentation of the manuscript. All mathematical arguments and results
were carefully verified by the authors, who take full responsibility
for the content of the paper.

\begingroup
\small
\bibliographystyle{plain}
\bibliography{references}
\endgroup

\par\bigskip
\begingroup
\footnotesize
\interlinepenalty10000
\indent\textsc{Einstein Institute of Mathematics, The Hebrew University of Jerusalem, Jerusalem, Israel}\par\nobreak
\indent\textit{Email addresses:} \nolinkurl{motty.gruda@mail.huji.ac.il},
\nolinkurl{yuval.peled@mail.huji.ac.il}\par
\endgroup

\end{document}